\documentclass[14pt]{article}%{book}%
\usepackage[english]{babel}
\usepackage[cp1251]{inputenc}
\usepackage{amsmath,latexsym}
\usepackage[psamsfonts]{amssymb}
\usepackage[dvips]{graphics}
\usepackage[mathcal]{euscript}
\begin{document}

\begin{center}\large
{\bf  FOLIATION THEORY AND IT'S APPLICATIONS } \\[4mm]
{\bf A.Ya.Narmanov,G.Kaypnazarova}\\[3mm]
{\bf National University of Uzbekistan}\\[3mm]
{\it narmanov@yandex.ru,gulita@rambler.ru}\\
\end{center}\vspace{3mm}

\begin{abstract}
Subject of present paper is the review of results of authors on
foliation  theory and  applications of foliation  theory  in
control systems. The paper consists of two parts. In the first
part the  results of authors on foliation theory are presented, in
the second part  the  results on applications of foliation theory
in the qualitative theory of control systems are given. In paper
everywhere smoothness of a class $C^\infty$ is considered.

\end{abstract}



\nopagebreak\footnotetext{Research supported by grant OT-F1-096 of
the Ministry of higher and secondary specialized education of
Republic of Uzbekistan.}
 \large
\begin{center}
\textbf{1. Topology of foliations}
\end{center}

The foliation  theory  is a branch of the geometry which has
arisen in the second half of the XX-th century on a joint of
ordinary differential equations and the differential topology.
Basic works on the foliation theory  belong to the French
mathematicians A. Haefliger \cite{HA1956},\cite {HA1958}, G.
Ehresman \cite{EH1961},\cite{EH1951}, G.Reeb
\cite{RE1952},\cite{RE1961}, H. Rosenberg \cite {Ro1966},\cite
{Ro1972},G. Lamoureux \cite{L1976}, \cite{L1984},R. Langevin
\cite{LAN1983},\cite{LP1993}. Important contribution to foliation
theory was made by known mathematicians - as well as I.Tamura
\cite{Ta1979}, R. Herman \cite{HE1960}
\cite{HER1960},\cite{HER1962} ,\cite{HE1962},T.Inaba [14],[15],
W.Turston [54],[55], P. Molino [26], P.Novikov [42], Ph. Tondeur
[56],[57],B. Reinhart [47]. Now the foliation theory is
intensively developed, has wide applications in various areas of
mathematics - such, as the optimal control theory, the theory of
dynamic polysystems. There are numerous researches on the
foliation theory.

The review of the last scientific works on  the foliation theory
and very big bibliography is presented in work of Ph. Tondeur
\cite {To1988}.

{\bf Definition-1.1} Let $(M,A)$ be a  smooth manifold of
dimension $n$, where $A $ is a $C^r$ atlas, $r\ge1$, $0<k<n$. A
family $F=\{L_\alpha:\alpha\in B\}$   of path-wise connected
subsets of $M$ is called $ k$-  dimensional $C^r-$ foliation of if
it satisfies to the following
three conditions:\\

$F_I :\ $ $ \bigcup \limits_{\alpha \in B} L_\alpha=M $  ;\\

$F_{II}:$ for  every $ \alpha,\beta \in B $ if  $ \alpha \neq
\beta $ then
$ L_\alpha \bigcap L_\beta = \emptyset$;\\

$F_{III}:$ For any point $p \in M $ there exists a local chart (
local coordinate system ) $ (U, \varphi) \in A, \ p \in U $ so
that if $ U \bigcap L_\alpha \neq \emptyset$ for some $\alpha \in
B $ the components of  $ \varphi (U \bigcap L_\alpha)$ are
following subsets of parallel affine planes
$$ \{ ( x_1, x_2, ..., x_n ) \in \varphi(U) : x_{k+1}= c_{k+1},
x_{k+2}= c_{k+2},...,x_{n}= c_{n} \}$$
where numbers $ c_{k+1}, c_{k+2},...,c_n $ are constant on
components (Figure-1,\cite {Ta1979}, p. 121).

The most simple examples of a  foliation  are given by integral
curves of a vector field and  by level surfaces  of differentiable
functions. If the $ X $ vector field without singular points  is
given on manifold $ M $ under the theorem of straightening of a
vector field (under the theorem of existence of the solution of
the differential equation) integral  curves generate
one-dimensional foliation on $M $.
\begin{center}
\includegraphics [16cm,8cm]{figure1.bmp} \\
                                 {Figure-1}
\end{center}

Let $M$ be a  smooth manifold of dimension $n$, $f:M\to R^1$ be a
differentiable function. Let $p_0\in M$,$ f(p_0)=c_0 $ and the
level set $L=\{p\in M:f(p)=c_0\}$ does not contain critical
points. Then the level set is a smooth submanifold of dimension
$n-1$. If we will assume that differentiable function has no
critical points, partition of $M$ into  level  surfaces of
function is a $n-1$- dimensional foliation (codimension one
foliation).Codimension one foliations generated by level surfaces
were studied in papers \cite {B2004},\cite {K2004},\cite {K2006},
\cite {KN2008},\cite {NBU2003}, \cite {NB2003},  \cite{NK2010},
\cite {Rei1959}, \cite {To1988}. The following theorem gives to us
a simple example of foliation.

\textbf{Theorem-1.1.} Let $f:M\to N$ be a differentiable mapping
of the maximum rank, where $M$ is a smooth manifold of dimension
$n$,$N$ is a smooth manifold of dimension $m$, where $n>m$. Then
for each point $q\in N$ a level set $L_q= \{p\in M: f(p)=f(q)\}$
is a manifold of dimension $n-m$ and partition of $M$ into
connection components of the manifolds $L_q$ is a $n-m$-
dimensional foliation.

Using the condition  3 of definition 1.1 it is easy to establish
that there is a differential structure on each leaf such that a
leaf is immersed $k$-dimensional submanifold of $M$, i.e the
canonical injection is a immersing map(a map of the maximum rank).
Thus on each leaf there are two topology: the topology $\tau_M$
induced from $M$ and it's own topology $\tau_F$ as a submanifold .
These two topologies are generally different. The topology $\tau
_F$ is stronger than topology $\tau _M$, i.e. each open subset
of $L_\alpha$ in topology $\tau_M$ is open in $\tau_F$.\\
A leaf $L_\alpha$ is called  compact if $(L_\alpha, \tau_F)$ is
compact topological space. It is obvious that the compact leaf is
a compact subset of  manifold $M$ . The leaf  is $L_\alpha$ called
as proper if the topology $\tau_F$ coincides with the topology
$\tau_M$ induced from $M$ .If these two these topology on
$L_\alpha$ do not coincide, the  leaf is called  a non- proper
leaf. It is easy to prove that the compact leaf is proper  leaf.
In work \cite{AN2004} the following assertion is proved which
takes place  for foliations with singularities too  which we will
discuss in the second part of this paper.

\textbf{Proposition.} If a  leaf  is a closed subset of $M$ then
it is a proper  leaf.

Let $L$ be a  leaf of $F$. Point $y\in M$  is called a limit point
of the leaf $L$ if there is a sequence of points $y_m$ from  $L$
which converges to $y$ in topology of manifold $M$  and does not
converge to this point in the  topology of the  leaf $L$ [36].

The set of all limit points of the  leaf $L$ we will denote by
$\Omega (L)$. It is easy to show that the limit set consists of
the whole leaves, i. e. if $y\in\Omega (L)$ then $L(y)\subset
\Omega (L)$, where $L(y)$ is a leaf containing $y$. Generally the
set $\Omega(L)$ can be empty or can coincide with all manifold. It
can already take place for trajectories of dynamic systems. For
example, if $L$ is closed the set $\Omega (L)$is empty, and in
case of an irrational winding of torus each trajectory everywhere
is dense and consequently its limit set coincides with all torus.
Studying of limit sets of leaves of foliation includes  studying
of limit sets of trajectories of the differential equations and it
is the important problem of the foliation theory. On these
subjects there are numerous
researches\cite{AN2004},\cite{In1977},\cite{M1948},\cite{N1981}-
\cite{N2000}. In work \cite{AN2004} following properties of a leaf
are proved which also takes place for foliation with singularities
too.

\textbf{Theorem-1.2}. (1). A leaf $L_0$  is proper leaf if and
only if   $L_0\bigcap \Omega(L_0)=0$;

(2). A leaf is $L_0$ is not proper leaf  if and only if
$\Omega(L_0)=\overline {L}_0 $ where $\overline{L}_0$ is the
closure of $L_0$  in manifold $M$.

For two leaves  $L_1$ and $L_2$ we will write in $ L_1\leq L_2$
only in a case when $L_1\subset\Omega(L_2)$. The inequality
$L_1<L_2$ means $L_1\leq L_2$ and $L_1\neq L_2$. The relation
$\leq$ on the set of leaves has been entered by the Japanese
mathematician T.Nisimori in the paper \cite{NISH1977}.

We will denote by  $(M/F,\leq)$ set of leaves with the entered
relation on it. It is obvious that the $\leq $ on $M/F$ reflective
and is transitive, but in many cases this relation is not
asymmetric, therefore generally the set $(M/F,\leq)$  is not
partially ordered.T.Nisimori was interested in the case where
$(M/F,\leq)$  is a  partially ordered set. Except that T.Nisimori
has entered concepts of depth of a leaf $L$  and depth of
foliation $F$ as follows: $ d(L)=\sup\{ k:$ there exist leaves $
L_1,L_2,...,L_k  $ such that $ L_1<L_2<...<L_k=L\}$, $
d(F)=\sup\{d(L): L\in M/F\}$.

A leaf $L$, being the closed subset, has the depth equal to one.
It is easy to construct one dimensional foliation of Euclid plane
with leaves of depth equal to two.

In work of \cite{NISH1977} Nishimori has proved the following
theorem which shows that for each  positive integer $k$ there
exists two-dimensional foliation with leaves of the depth equal to
$k$.

\textbf{Theorem-1.3.} Let $S_2$ be a closed surface of a genus
2.For all positive integer $k$ there is a  codimension one
foliation $F$ on $M=S_2\times [0,1]$ satisfying the following
conditions (1),(2)and(3).

(1) All leaves of $F$ are proper and  transverse to
${x}\times[0,1]$ for all $ x\in S_2 $. $S_2\times{0}$ and $
S_2\times{1}$ are compact leaves.

(2) $d(F)=k$.

(3) All holonomy groups of $F$ are abelian.

The following theorem is proved in paper \cite{AN2004} shows that
there exists one dimensional analytical foliation generated by
integral curves of analytical vector field which have leaves of
depth equal to 1,2 and 3.

\textbf{Theorem-1.4}. Let $S^k$ be a $k$ dimensional sphere.On the
manifold $ S^2\times S^1 $ there exists  an analytical vector
field without  singular points and with three  pairwise different
integral curves $\alpha,\beta,\gamma$ such that
$\alpha\subset\Omega(\beta)$,$\beta\subset\Omega(\gamma)$,where
$\alpha$ is a closed trajectory, $\Omega(\beta) $ consists of only
closed trajectories,$\Omega(\gamma)$ consists of only  the
trajectories of depth equal to two.

This vector field generates one-dimensional foliation of the depth
equal to 3.

\textbf{Remark.}The example of not analytical dynamic system of a
class $C^{\infty}$ was constructed in the paper \cite{M1948}  for
which there is an infinite chain of not closed trajectories $ L_i
$ such that  each trajectory  $ L_{i+1} $ is in the limit set of $
L_i $.

In the paper \cite{NISH1977} for codimension one foliation  the
following theorem is proved:

\textbf{Theorem-1.5.} (Nishimori).If $ d(F)<\infty$  or all leaves
of  foliation $F$ are proper,then the  set  $(M/F,\leq )$ is
partially ordered.

Nishimori, studying property codimension one foliation in the case
when the set $(M/F,\leq)$ is a partially ordered, has delivered
following questions which are of interest  for  foliation with
singularities too \cite{NISH1977}: \\
1. Are all leaves of foliation $F$ proper under the assumption
that the set$(M/F,\leq)$
is partially ordered? \\
2. Is a leaf  $L$ proper under the assumption that $dL<\infty$?

A.Narmanov studied the relation $\leq $ for foliation with
singularities in the paper \cite {N1981}. In particular, he proved
the following theorems which solves problems 1, 2 delivered by
Nishimori.

\textbf{Theorem-1.5.} Let $M/F$ be the  set of leaves of
foliation $F$ with singularities .Then the set  $(M/F,\leq)$ is a
partially ordered if and only if all leaves are proper.

\textbf{Theorem-1.6.}If  the depth of a leaf is finite, then it is
proper leaf.

It is known that the limit leaf  of compact leaves codimension one
foliation on compact manifold is a compact leaf and the limit set
of each leaf contains finite number of compact leaves.

The following theorems are generalizations of these facts for
leaves with finite depth \cite {N1983}.

\textbf{Theorem-1.7.}Let $F$ be a transversely oriented
codimension one foliation on compact manifold $M$, $L_i$- a leaf
of  foliation  $F$, and $x_i\longrightarrow x$, where $x_i\in
L_i$. If $dL_i\leq k $ for each $i$ then $ dL(x)\leq k $.

\textbf{Theorem-1.8.} Let $F$ is a transversely oriented
codimension one foliation on compact manifold $M$,$L^{0}$ - some
leaf of foliation $F$. Then for each $k\geq1$ the set
$C_k=\{L:L<L^{0},dL=k\}$ either is empty, or consists of finite
number of leaves.

Let's remind that transversally orientability of  $F$ means that
there exists smooth non-degenerated vector field $X$ on $M$, which
is transversal to leaves of  foliation $ F$.

Let $ x \in M $, $L(x) $ is a leaf foliation,containing the point
$x$, $T_x$ is a manifold dimension of $n-k$ transversal to $L(x)$
such that $ T_x \bigcap L(x)=x $.To each the closed continuous
curve in  $ L(x) $ beginning and the ending at the point $ x\in M$
corresponds a local diffeomorphism $ g $ of the manifold $ T_x $,
given in some neighborhood of the  point $ x $  in $ T _x $ such
that $ g(x)=x $. The set of such diffeomorfisms forms the
pseudogroup $ \Gamma_x(L) $ of the leaf $ L $ at the point $x$,
and germs of these diffeomorphisms form holonomy  group $ H $ of
the leaf $ L(x)$.For different points from $L$ corresponding
holonomy groups  are isomorphic\cite {Ta1979}.\\
The important results in foliation theory  are received by G.
Reeb. One of his theorems is called as the theorem of local
stability which can be formulated as follows.

\textbf{Theorem-1.9.}\cite {Ta1979} Let $ L_0 $ a compact leaf
foliation $ F $ with finite holonomy group. Then there is an open
saturated set  $ V $ which contains $ L_0 $ and consists of
compact leaves.

Let's notice that a saturated set $ S\subset M $ on a foliated
manifold is a subset which is the union of leaves.

In 1976 in Rio de Janeiro at the international conference the
attention to the question on possibility of the proof of theorems
on local stability for noncompact leaves \cite {Sch1978} has been
brought. In 1977 the Japanese mathematician T.Inaba has
constructed a counterexample which shows that if codimension of
foliation is not equal to one G.Reeb's theorem cannot be
generalized for noncompact leaves \cite{In1977}.

Let's bring  the theorem on a neighborhood of a leaf with finite
depth which is generalization of the theorem of J.Reeb on local
stability for transversely oriented codimension one foliation.

Let $F$ transversely oriented, $X$ a smooth vector field on $M$,
transversely to leaves of $F$. Let $x\in M,t\rightarrow X^{t}(x)$,
- the integral curve of a vector field $X$ passing through the
point $x$ at $t=0$. Let's put $T_x=\{X^{t}(x):-a<t<a \}$. In
further will write $T_x\approx (-a,a)$ and as usually, to replace
subsets of $ T_x $ by  their images from $(-a,a)$. The point is
$y\in T_x$ called as a motionless point of pseudo-group
$\Gamma=\Gamma _x(L)$, if $g(y)=y$ for each $g\in \Gamma$,
advanced in a point $y\in T_x$. If there exists $\varepsilon > 0$
such that each point from $( -\varepsilon,\varepsilon)$ is  a
motionless point of pseudogroup $\Gamma$ we will say that the
pseudo-group is $ \Gamma$ trivial.

Let $F$ be a codimension one foliation, $L$ be a some leaf of $F$
with finite depth ,$\rho$ - distance function defined by some
fixed riemannian metric on $M$.

Let's enter set $ U_r=\{y\in M:\rho(y,L)<r\},r>0 $, where $
\rho(y,L) $- distance from the point  $ y $ to the leaf $ L $.

\textbf{Theorem-1.10}(\cite {ND1983}). Let $ F $ be a transversely
oriented codimension one foliation on compact manifold $M$. If the
holonomy pseudogroup $ \Gamma $ the leaf $L$ is trivial,then for
each $r>0$ there is a invariant open set $V$ containing $L$ and
consisting of leaves diffeomorphic
to $ L $ which satisfies to following conditions: \\
1)$ V\subset U_r $;\\
2) $ dL_\alpha=dL $ for each leaf  $ L_\alpha\subset V $.

One more generalization of G. Reeb theorem   for a noncompact leaf
is resulted below.For this purpose we will bring some
definitions.\\
Let $ M $ be smooth connected complete  riemannian manifold of
dimension$ n $ with riemannian metric $ g $, $ F $ -smooth
foliation of  dimensions $ k $
on $ M $.\\
Let's  denote through $ L(p)$ a leaf of $ F $ passing through a
point $ p $, $ F(p) $- tangent space of leaf  at the  point $p$,
$H(p)$ - orthogonal complementary of  $ F(p)$ in  $ T_p M $, $p\in
M $. There are two subbundle (smooth
distributions),$TF=\{F(p):p\in M\} $,$ H=\{H(p):p\in M\}$of
tangent bundle $ TM $ such ,that   $ TM=TF\oplus H $ where is $ F
$ orthogonal addition $ TF $.\\
Piecewise smooth curve  $ \gamma
:[0,1]\rightarrow M $ we name horizontal, if  $\dot{\gamma}(t)\in
H(\gamma(t)) $ for each $t\in [0,1] $.Piecewise smooth curve which
lies in a leaf foliation  $F$ is called as vertical.

Let $I =[0,1] $, $ \nu :I\rightarrow M $  a vertical curve,$ h:
I\rightarrow M $ a horizontal curve and $ h(0)=\nu (0)$.Piecewise
smooth mapping $ P:I\times I\rightarrow M $ is called as
vertical-horizontal homotopy for  pair $v,h $ if   $ t\rightarrow
P(t,s) $ is a vertical curve for each $ s\in I $, $ s\rightarrow
P(t,s) $ is a horizontal curve for each $ t\in I $,and  $
P(t,0)=\nu(t) $ for $ t\in I $, $ P(0,s)=h(s) $ for $ s\in I $. If
for each pair of vertical and horizontal curves  $ \nu, h
:I\rightarrow M $ with  $ h(0)=\nu (0) $ there  exists
corresponding vertical-horizontal homotopy  $ P $ ,we say that
distribution  $ H $ is  Ehresman connection  for $ F $
\cite{BHE1984}.

Let $ L_0 $ a leaf of codimension one foliation $ F $, $ U_r=\{
x\in M: \rho (x;L_0)< r\}$, where $ \rho (x,L_0) $- distance from
the point$ x $ to a leaf $ L_0 $. We will assume that there is
such number $ r_0 >0 $ that  for each horizontal curve $ h:
[0,1]\rightarrow U_{r_0}$ and for each vertical curve $ \nu
:[0,1]\rightarrow L_0 $ such that $ h(0)=v(0) $ there exists
vertically-horizontal homotopy for  pair $(\nu,h ) $. At this
assumption  we formulate generalization of the theorem of J.Reeb
\cite{N1995}.

\textbf{Theorem-1.11.} Let $ F $ a transversely oriented
codimension one foliation , $ L_0 $ be a relatively compact
proper leaf leaf with finitely generated fundamental group. Then
if holonomy group  of the leaf $ L_0 $  is trivial then  for each
$ r>0 $ there is an saturated set $ V $ such that  $ L_0 \subset V
\subset U_r $ and restriction of  $ F $ on  $ V $ is a fibrarion
over  $ R^{1}$ with the leaf $ L_0 $.

From the geometrical point of view, the important classes of
foliation are total geodesic  and riemannian foliations. Foliation
$F$ on  riemannian manifold $ M $ is called  total geodesic  if
each leaf of foliation $F$ is a total geodesic submanifold,i.e
every geodesic  tangent to a leaf   foliation $ F $ at one point,
remains on this leaf. The geometry of total geodesic   foliations
is studied in works \cite{HE1962}, \cite{JW1979},
\cite{N1999},\cite{BH1984}.

Foliation $F $ on a  riemannian manifold $ M $ is called
riemannian if each geodesic, orthogonal at some point to a leaf of
foliation $F $ , remains orthogonal at all points to leaves of $F
$ \cite{Rei1959}. Riemannian foliation without singularities for
the first time have been entered and studied by Reinhart in work
\cite{Rei1959}. This class foliation naturally arising at studying
of bundles and level surfaces. Riemannian foliation are studied by
many mathematicians, in particular,in works of R.Herman
\cite{HE1960},\cite{HER1960},\cite{HER1962},P.Molino
\cite{MO1988},A Morgan \cite{MOR1976}, Ph.Tondeur \cite{To1988} .
The most simple examples of Riemannian foliation are partition of
$ R^{n} $ into parallel planes or into concentric hyperspheres.
Riemannian foliation with singularities have been entered and
studied in works of P.Molino \cite{MO1988}, they also were studied
by A.Narmanov in the papers \cite{N1999},\cite {N1996}.

Let $ M $ smooth connected complete  riemannian reducible
manifold. Then on $ M $ there are two parallel foliations,
mutually additional on orthogonality \cite{KN1981}. If $ M $
simple connected manifold then the de Rham theorem takes place
which asserts  that $ M $ is isometric to direct product of any
two leaves from different foliations \cite{KN1981}. In this case
both foliations are riemannian and total geodesic simultaneously.
Below it is presented results of authors on geometry of riemannian
and totally geodesic foliations.

Assume that $ F $ is a riemannian foliation with respect to
riemannian  metric  $ g $ on $ M $. Let  $ \pi_1 : TM \rightarrow
TF $, $ \pi_2 :TM \rightarrow H $ be orthogonal  projections, $
V(M) $ ,$ V(F)$, $ V( H) $ be the  sets of smooth sections of
bundles  $ TM $ ,$ TF $ and $H$ accordingly. If $X\in V(F) $
($X\in V(H) $)  $ X$ is a called a vertical (horizontal) field.

Now we will assume that each leaf of  $ F $ is a total geodesic
submanifolds of $M$. It is equivalent to that, $\nabla _X Y  \in
V(F)$ for all $ X,Y \in V(F) $ \cite {Rei1959}, where $\nabla $-
Levi-Civita connection. In this case $ F $ is a riemannian
foliation with total geodesic leaves. Then on bundles $ TF $ and $
H $ are given metric connections $\nabla^{1}$ and $ \nabla^{2}$ as
follows.If $X\in V(F)$,$ Y\in V(H) $, $ Z\in V(M),$ we will put
$$\nabla ^1_Z=\pi _1 (\nabla _Z(X)),\   \nabla
^2_Z(Y)=\pi _2[Z_1,\widetilde{Y}]+\pi _2[\nabla _{Z_{2}},Y],$$
where $ Z=Z_{1}\oplus Z_{2}$,$Z_{1}\in V(F),Z_{2}\in
V(H),\widetilde{Y}\in V(M) ,\pi _{2}\widetilde{Y}=Y$. Here
$[Z_{1},\widetilde{Y}]$ - Lie bracket of vector fields $Z_{1}$ and
$ \widetilde{Y}$ . $ \nabla _{2}$ is a metric connection  only in
case when $ F $ is riemanian. Owing to that $F$ is a totally
geodesic foliation, connection  $ \nabla^{1}$ also is metric
\cite{N1999}, \cite{Rei1959}.

Let $ p\in M $, $ S(p) $ be the  set of points $M$ ,which can be
connected by horizontal curves with $p$. Owing to that foliation
$F$ is total geodesic, for each the $ p\in M $  set $ S(p) $ has
topology and differentiable structure, in relation to which $ S(P)
$ is a immersed submanifold of $M$ \cite{BH1984}.It is easy to
prove the following assertion.

\textbf{Lemma-1.1.} $ dimS(p)\geq k $ for every point $p\in M$.

Owing to that the manifold $M$ is complete, and considered
foliation is a riemannian, the distribution $H$ is a Ehresmann
connection for $F $ \cite{BHE1984}. Therefore for each piece-wise
smooth curve $\gamma :I\to M$ there exists unique
vertical-horizontal homotopy $P_\gamma :I\times I\to M$ such that
$\gamma (t)=P_\gamma (t,t)$ for $t\in I$. Let $P:I\times I\to M$
be a vertical-horizontal homotopy. We will denote by $D_t
(P(t,s))$ the tangent vector of the curve $t\to P(t,s)$ at the
point $P(t,s)$, by $D_s (P(t,s))$ the tangent vector of the curve
$s\to P(t,s)$ at the point $P(t,s)$.

\textbf{Lemma-1.2.}  Let $X(t,s)=D_s (P(t,s))$, $Y(t,s)=D_t
(P(t,s))$ for $(t,s)\in I\times I$. Then $ \nabla ^1_X Y=0 $ and $
\nabla ^2_Y X=0 $.

In the proof of the de Rham theorem projections of a curve
$\gamma$ in $ L(p_0) $  and in $ S(p_0) $ are defined with
connection $ \nabla $  and it is shown that these projections
coincide with curves $ P_\gamma (\cdot ,0 ) :I \rightarrow M $ and
$ P_\gamma (0 ,\cdot ) :I \rightarrow M $ accordingly. In this
case it is used that distribution $H $ is  complete integrable
\cite{KN1981}. In the paper \cite{N1999} the similar fact is
proved for projections of a curve $\gamma$ without the assumption
that  $H $ is  complete integrable by means of metric connection
$\widetilde\nabla$ which is entered below. Metric connection
$\widetilde\nabla$ is defined as follows:
$$\widetilde{\nabla } _Z X=\nabla ^1_Z X_1 +\nabla ^2_Z X_2,$$
where $X,Z\in V(M)$, $X_i=\pi_i(X)$, $i=1,2$.It is not difficult
to check up that distributions $TF$ and $H$ are parallel with
respect to $\widetilde{\nabla }$.

Let $ \gamma :I\to M $ be  a smooth curve,$ \gamma(0)=p_0 $ and
$\gamma(1)=p$, $ C:I \to T_{p_0} M $ is a  development of the
curve $ \gamma $ in  $ T_{p_0} M $ defined by connection
$\widetilde\nabla$. (See definition of development in
[\cite{KN1981},p.129]. (Here for convenience tangent vector space
$ T_{p_0}M $ is identified with affine tangent space at the point
$p_0.)$

Let $ C(t)=(A(t),B(t)) $ where $ A(t)\in F(p_0)$,$ B(t)\in H(p_0)$
for $ t\in I $. As $M$ is a complete riemannian manifold,
$\widetilde\nabla$- metric connection, there are smooth curves  $
\gamma_1,\gamma_2 :I\to M $, which are developed on curves $ t\to
A(t)$ and $ t\to B(t) $ accordingly (\cite{KN1981}, p.167, the
theorem 4.1). According to the proposition 4.1 in
(\cite{KN1981},p.129) $\gamma_i$ is a such curve that result of
parallel transport $\dot\gamma_i$ to the point $p_0$ along
$\gamma_i^{-1}$ defined by connection $\widetilde\nabla $
coincides with the result of parallel transport of
$\pi_i(\dot\gamma(t))$ to the point $p_0$ along $\gamma^{-1}$
defined by connection $\widetilde\nabla $ too.That is why
$\gamma_1$ is a vertical curve, $\gamma_2$ is a horizontal curve.
Curves $\gamma_1$, $\gamma_2$ are called projections of the curve
$\gamma$  in $L(p_0)$ and in $S(p_0)$accordingly.

The following theorems are proved in the work \cite{N1999}.

\textbf{Theorem-1.12.} The projection of  curve $\gamma:I\to V$ in
$L(q_0)$ (in $ S(q_0)$) is a curve $ P(\cdot,0):I\to L(q_0)$
(accordingly $P(0,\cdot):I\to S(q_0)$ ).

The following theorem shows that if distribution $ H $ is complete
integrable if and only if connection  $\widetilde{\nabla }$
coincides with Levi-Civita  connection $ \nabla $ .

\textbf{Theorem-1.13.} Following assertions are equivalent.\\
1.Distribution  $H$  is complete integrable (i.e  $dimS(p)=n-k$ for each $p\in M$).\\
2. $\widetilde{\nabla }$ is connection without torsion (i.e $
\widetilde{\nabla }=\nabla $ ).

\textbf{Remark }. As shows known Hopf fibration of  on
three-dimensional sphere,the  distribution $ H $ it is not always
complete integrable. In a case when $ H $ is complete integrable
de Rham theorem  takes place: if $ M $ is simple connected,it is
isometric to product $ L(p)\times S(p) $ for each $ p\in M $. In
this case $S(p)$ is a leaf of foliation $F ^\bot$ generated by
distribution $H$. Projections of any point $ p\in M $ in $L(p_0)$
and  in $S(p_0)$ are defined  as follows. Let $\gamma:I\to M$ is a
smooth curve, $\gamma(0)=p_0, \ \gamma(1)=p ,\ \nu $, $h$ are
projections of  $\gamma$ in $L(p_0)$ and in $S(p_0)$. The points
$\nu(1)$, $h(1)$ are called as projections of $p$ in $L(p_0)$ and
in $S(p_0)$ accordingly. Owing to that the distribution $H$ is
completely integrable, the projection of $p$ depends only on the
homotopy class of the curve $\gamma$.That is why  when $M$ is
simple connected,the mapping $f:p \to(p_1,p(_2)$ is correctly
defined. Under theorems 1.12 and 1.13 mapping $f$ is a isometric
immersing. Since $dimM=dim\{L(p_0)\times S(p_0)\}$ the  mapping
$f$ is covering mapping, hence, it is an isometry
(\cite{KN1981},p.134).

In the known monograph \cite{To1988} Ph. Tondeur studied
foliation, generated by level surfaces of   functions of a certain
class. He considered function $f:M \to R ^1$ without critical
points on Riemannian manifold $M$ for  which length of a gradient
is constant on each level surface. For such functions he has
proven that foliation generated by level surfaces  of such
function, is a Riemannian foliation . Authors of the present
article studied geometry of foliation generated by  level surfaces
of the functions considered in the monograph of professor Ph.
Tondeur without the assumption of absence of critical points.

\textbf {Definition -1.2.}Let $M$ be a smooth manifold of
dimension $n$.Function $f:M \to R ^1$  of the  class  $C^2
(M,R^1)$ for which length of a gradient is constant on  connection
components of level sets   is called a metric function.

Let $f:R^n \to R ^1$ be a metric function. We will consider system
of the differential equations

\begin{center}
$\dot{x}=gradf(x)$  \vspace{1mm}                         (1.1)
\end{center}\vspace{5mm}
As, the right part of system (1.1) is differentiable, for each
point $x_0\in M$ there is a unique solution of system (1.1) with
the initial condition $x(0)=x_0$. The trajectory of system (1.1)
is called gradient line of the  function $f$.

\textbf{Theorem-1.14}\cite{K2004}.Curvature of each gradient line
of metric function is equal to zero.

In the papers \cite{K2004},\cite{KN2008} topology of level
surfaces are studied under the assumption that each of a
connection components of the set of critical points of metric
function is either a point, or is a regular surface and every
component is isolated from others.

The following theorem gives complete classification of foliations
generated by level surfaces of metric function \cite{KN2008}.

\textbf{Theorem-1.15.} Let $f:M \to R ^1$  be a metric function
given in $R^n$. Then level  surfaces of  function $f$  form
foliation which has one of following $n$ types :\\
1) Foliation $F$ consists of parallel hyperplanes;\\
2) Foliation $F$ consists of concentric hyperspheres and the point (the center of hyperspheres); \\
3) Foliation $F$ consists of concentric cylinders of the  kind
$S^{n-k-1}\times R^k$ and the singular leaf $R^k$ (which arises at
degeneration of spheres to a point), where $k$- the minimum of
dimensions of critical level surfaces,$1\leq k \leq{n-2}$.

At the proof of the theorem-1.1 4 the following theorem is used
which also is proved in work \cite{KN2008}.

\textbf{Theorem-1.16.} Let $L$ be  a regular surface of dimension
$r$ which is the closed subset $R^n$, where $1\leq r\leq {n-1}$.
If the normal planes passing through various points $L$ are not
crossed,then $L$ is $r$- a dimensional plane.

This theorem represents independent interest for the course of
differential geometry.

The following theorem shows internal link between geometry of
riemannian manifold and property of metric function which is given
on it.

\textbf{Theorem-1.17} \cite{NK2010} Let $(M,g)$ be a smooth
riemannian manifold of dimension $n$,  $f:M \to R ^1$ metric
function. Then each gradient line of the  function  $f$ is a
geodesic line  of the riemannian manifold $M$.

For the metric functions given on a riemannian manifold , it is
difficult to get classification theorems, as in a case $M=R^n$.
Here one much depends on topology of  riemannian manifold on which
function is given. For example, in Euclid  case if metric function
has no critical points, then as shown \cite{KN2008} all level
surfaces are hyperplanes. It is easy to construct metric function
without critical points on the two-dimensional cylinder ( with the
metric induced from Euclid structure of  $R^3$)  level lines of
which are circles (compact sets).

The following theorem, is the classification theorem for level
lines   of the metric function given on two-dimensional riemannian
manifold.

\textbf{Theorem-1.18}\cite{NK2010} Let $M$ be two-dimensional
riemannian manifold,  $f:M \to R ^1$ be a  metric function without
critical points. Then all level lines  are homeoomorphic to
circle, or all level lines  are homeoomorphic to a straight line.

The following theorem shows that if the metric function is given
on complete simple connected riemannian manifold and does not have
critical points, then it has no compact level surfaces.

\textbf{Theorem-1.19} \cite{NK2010}  Let $M$ be a smooth complete
simple connected riemannian manifold, $f:M \to R ^1$ be a metric
function without critical points. Then level surfaces are mutually
diffeomorphic noncompact submanifolds of $M$.

\begin{center}
\textbf{2. Applications of foliation theory in control systems}\\
\end{center}

The last years methods and results of foliation theory  began to
be used widely in the qualitative theory of optimal control.It was
promoted by works of the American mathematician of G.Sussmann
\cite{Su1973} and the English mathematician P.Stefan \cite{St1974}
which have shown that a orbit of family of smooth vector fields is
a smooth  immersed submanifold. Besides, they have shown that if
dimensions of orbits are the same, partition of phase space into
orbits is a foliation. Papers \cite{N1981}, \cite{ND1983},
\cite{N1995}, \cite{N1997}, \cite{GB1982}, \cite{N2000},
\cite{N1996}, \cite{NPI1985}, \cite{NP1985} are devoted to
applications of foliation theory in control theory.

Let $D$ be a family of the smooth vector fields  on $M$, and $X\in
D$. Then for a point $x\in M$ by $X^t(x)$ we will denote the
integral curve of a vector field $X$ passing through the  point
$x$ at $t=0$. Mapping $t\to X^t(x )$ is defined in some domain
$I(x)$ which generally depends not only on a field $X$, but also
from the point $x$. Further everywhere in  formulas kind of
$X^t(x)$ we will consider that $t\in I(x)$. The orbit
$L(x)$passing through a point $x$ of the family of  the vector
fields $D$ is defined as a set of  points $y$ from $M$ for which
there are real numbers $t_1,t_2,\ldots,t_k$ and vector fields
$X_{i_1},X_{i_2},\ldots X_{i_k}$ from (where $k$ is  a natural
number) such that $y=X^{t_k}_{i_k}(X^{t_{k-1}}_{i_{k-1}}(\ldots
(X^{t_1}_{i_1})\ldots ))$.

Now we will bring a definition of foliation with singularities
\cite{St1974}. A subset $L$ of manifold  $M$ is called as a leaf
if

1) there is a differential structure $ \sigma $ on $ L $ such that
$(L,\sigma)$ is a connected  $k$- dimensional immersed submanifold of  $M$. \\
2) for locally connected topological space  $ N $ and for
continuous mapping $ f:N\to M $ such that $ f(N)\subset L$ the
mapping  $f:N\to (L,\sigma)$ is continuous.

Partition $F$ of manifold $M$ into leaves is called  smooth (of
the  class $C^r$ ) foliation with singularities if following
conditions are satisfied:

1) For each point  $x\in M$ there exists a local  $C^r$- chart
$(\psi,U)$ containing the point $x$ such that $\psi(U) =V_1\times
V_2$ where $V_1$ is a  neighborhood of origin in $R^k$, $V_2$- a
neighborhood of origin in $R^{n-k}$,$k$- dimension of the leaf
containing the point $x$;

2) $\psi (x)=(0,0)$ ;

3)For each leaf $L$ such that  $L\cap U\neq \emptyset$ it takes
place equality $L\cap U=\psi^{-1}(V_1\times l)$ where $l=\{\nu \in
V_2:\psi ^{-1}(0,\nu)\in L\}$.

By definition 1.1  each regular foliation  is a foliation  in
sense of the above-stated definition. In this case  every
connection component  of the set $l $ is a point. If dimensions of
leaves of a foliation with singularities are the same as noted
above, it is a foliation in sense of the definition 1.1. Thus, the
conception of foliation with singularities is a generalization of
classical notion of a foliation (now which  is called as regular
foliation). In the literature instead of "foliation with
singularities" the term "singular foliation"  is used also
\cite{MO1988}.To the studying of a foliation with singularities
are devoted papers \cite{AN2004},\cite{MO1988}, \cite{N1995},
\cite{N1997}, \cite{St1974}.

Now we will consider some applications of the foliation theory in
problems of the qualitative theory of control systems.

Let's consider a control system
\begin{center}
   $ \dot{x}=f(x,u)$,$ x\in M $, $ u\in U\subset R^{m}$                           $ (2.1)$
\end{center}
where $ M $ is a  smooth (class$ C^{\infty}$)  connected manifold
of  dimension $ n $ with some riemannian metric $ g $, $ U $ is a
compact set, for each the $ u\in U $ vector field $ f(x,u )$ is a
field of class $ C^{\infty}$, and mapping  $ f:M\times U
\longrightarrow TM $ ,where $ TM $  is the tangent bundle of   $ M
$, is continuously differentiable. It means that there is such
open set $ V $ such that $ U\subset V \subset R ^{m} $, and
continuously differentiable mapping $ \bar{f} : M\times
V\longrightarrow TM $, restriction of which on $ M\times U $
coincides with $ f(x,u )$. Admissible controls are defined as
piecewise-constant functions $ u:[0,T]\longrightarrow U $ ,
 where  $ 0<T<\infty $. Thus, the trajectories of system (2.1)
 corresponding to admissible controls, represent piecewise smooth mapping $ x:[0,T]\rightarrow M $.

The purpose of control is  a bring of the system to some fixed
(target) point $ \eta \in M $. We will say that the point
$x_{0}\in M$ is controllable from a point $\eta$ in time $T>0$, if
there is such trajectory of $x:[0,T]\rightarrow M$ of system (2.1)
that $x(0)=x_{0}, x(T)=\eta $. Let's denote by $G_{\eta}(<T)$ a
set of points of $ M $ which are controllable from a point $ \eta
$ for time, smaller than $ T $.We assume that $ \eta \in G{\eta}
(<T) $  for each  $ T> 0 $. The set of all points $ M $,which are
controllable from a point $ \eta $, is called as set  of
controllability with a target point  $ \eta $ and is denoted by $
G_{\eta}$. We will denote by  $ T=T_{\eta}(x) $ function of
Bellman given on set $G_{\eta}$ for the  optimal time problem. It
is known that a set of smooth vector fields on a smooth manifold
can be transformed into Lie algebra in which as product of vector
fields $ X $ and $ Y $ serves their Lie bracket $[ X,Y] $.

Let's denote by $ D $ set of vector fields $\{f(\cdot,u):u\in U\}
$,  by  $ A(D) $ minimal Lie subalgebra, containing $D$, by
$A_{x}(D)$ the subspace tangent spaces at a point $ x\in M $,
consisting of all vectors $\{X(x): X\in A(D)\} $. If we will
denote by $ L(\eta) $ a orbit of family $ D $ containing the point
$ \eta $, then it  follows from definition of the orbit  that  $
G_{\eta} \subset L(\eta ) $ for all $ \eta\in M $. The following
assertion \cite{NPI1985} takes place.

\textbf {Theorem-2.1}. If $ dimA_{\eta}(D)=dimL(\eta) $ then $ int
G_{\eta} \neq \varnothing $ in topology of $ L(\eta) $ .

Now let us give following definitions.

\textbf {Definition-2.1}. We will say that the system (2.1) is
completely controllable on $ L(\eta_{0}) $, if for all $ \eta\in
L(\eta_{0})$ it takes place equality $ G_{\eta}=L(\eta_{0})$ .

\textbf {Definition- 2.2}. The system (2.1) is called  normally
-locally controllable  (or, more shortly, $N$- locally
controllable) near  a point $ \eta $ if for any $T>0 $ there is a
neighborhood  $ V $ of the point $ \eta $ in $  L(\eta ) $ such
that each point from  $ V $ is controllable from  $ \eta $ in
time, smaller $ T $.

If the system is $ N $- locally controllable near  each point of $
L(\eta)$  we will say that it is  $ N $ - locally controllable  on
$ L(\eta)$ (see \cite{Su1973}).

\textbf {Definition - 2.3.} We will say that the system (2.1) is
completely ($ N $ -locally) controllable on invariant set $ S $ if
it is completely  controllable ( $ N $ - locally controllable) on
each a leaf of $ S $ .

Assume that  $ dimA_{x}(D)=k $ for every $ x\in M $, where
$0<k<n$,  $ A_{x} (D)=\{X(x):X\in A(D)\}$. In this case splitting
of $ F $ manifold $ M $ into  orbits family $ D$ is  a
$k$-dimensional foliation, i.e.orbits are $ k $ dimensional
submanifolds of  $M$.

 Let's consider the following question: if
the system (2.1) has property of complete controllability on one
fixed leaf of the foliation $ F $, under what conditions the
system (2.1) has this property on leaves close to a given leaf?

This question closely related with problems of the qualitative
theory of foliations on local stability of a leaf in sense of
J.Reeb(see \cite{Ta1979}). In the paper \cite{GB1982} the answer
is given to this question in the case when the leaf $ L_{0} $ of $
F $ in neighborhood of which the system (2.1) is studied,is
compact set. In this case conditions of the theorem of J.Reeb on
local stability is required. Namely, the following theorem is
proved.

\textbf{ Theorem-2.2}. Let  $ L_{0}$ be a compact leaf of $ F $
with finite holonomy group.If the system (2.1) is complete
controllable on $ L_{0} $ then it is is complete controllable
 on leaves,close to $ L_{0} $.

Thus, existence of such saturated(invariant) neighborhood $ V $ of
a leaf  $L_0$  gives the sufficient condition for stability of the
complete controllable system (2.1) on $ L_{0} $  , when $ L_{0} $
is a compact leaf. As Example 3 in \cite {GB1982} shows, complete
controllability on close leaves does not follow from the fact that
the system (2.1) is  complete controlled on a noncompact proper
leaf $L_0$ having a neighborhood $V$ described in theorem 2.1.
Therefore, in a case when $L_0$ is a noncompact leaf, we need
additional conditions that guarantee stability of the complete
controlled system (2.1) on $L_0$ . The theorem 1.11 gives the
possibility to get sufficient conditions for stability of the
complete controllable system (2.1) on $ L_{0} $ , when $ L_{0} $
is  a noncompact leaf.

\textbf {Theorem-2.3} \cite{NK2010} Let $ F $ be a transversely
oriented codimension one foliation, $ L_{0} $ be a relatively
compact proper leaf  with finitely generated fundamental group and
with trivial  holonomy group . Then if the system (2.1) is  $ N
$-locally controllable on $\overline {L}_0 $ (the closure in $ M
$) then there is an open saturated neighborhood $ V $   of the
leaf $ L_{0} $ such that on each leaf from $ V $ the system (2.1)
is complete controllable.

\textbf {Remark 1}. The closure of each leaf is an invariant set
(see (\cite{Ta1979},Theorem 4.9).

\textbf {Remark 2}. If the system (2.1) is $N$-locally controlled
on a leaf $L$ then it is complete controlled on $L$ (see
\cite{N2000}).

Let now  $ dimA_{x}(D)=k $ for every $x\in M$ where $0<k<n$, $F$
is a riemannian foliation with respect  to riemannian metric $g$.
We will remind that foliation $ F $ is called riemannian, if each
geodesic orthogonal at some point to a leaf foliation $ F $ ,
remains orthogonal to leaves at all points.

\textbf {Theorem-2.4} \cite{N2000} Let $(M,g)$ be a complete
riemannian manifold, $ L_{0} $ be a relatively compact proper leaf
of $ F $ .Then if the system (2.1) $ N $ -is locally controllable
on $\overline{L}_0 $ (on the closure of $ L_{0} $ in $M$)then
there is an invariant neighborhood $ V $ of the leaf $ L_{0} $
such that system (1) is complete controllable on each leaf from $
V $.

In the paper \cite{MO1988} it is given the necessary and
sufficient condition for singular foliation $ F $ to be riemannian
foliation. This condition deals with vector fields from $ A(D) $
and riemannian  metric $ g $.

Let now  $F$ be $k $-dimensional riemannian foliation with
respect  to riemannian metric $ g $, where $0<k<n$. It is possible
to present each vector field  $ X\in V(M) $ as $ X=X_{P}+X_{H} $
where $ X_{P},X_{H}$ orthogonal projections  of $X $ on $ P $ and$
H $ accordingly. If$ X_{H}=0 $, then  $ X\in V(F) $  and  the
vector field $ X$is called tangent vector field , if $ X_{P}=0 $
then $ X\in V(H)$ and the vector field $X$ is called horizontal
field. For vector fields $X,Y$ we will consider the bilinear
symmetric form $g_T(X,Y)=g(X_H,Y_H)$ on  $V(M)$, kernel of which
coincides with $V(F)$. We will study properties of this form. By
the definition of foliation for each point $p\in M$ there is a
neighborhood $U$ of the point $p$ and local system of coordinates
$x^1,x^2,\ldots x^k,y^1,y^2,\ldots y^{n-k}$ on $U$ such that
$\frac{\partial }{\partial x^1},\frac{\partial }{\partial
x^2},\ldots ,\frac{\partial }{\partial x^k}$ form basis of
sections of $ TF| _U $.The basis $\nu_{k+1},\nu_{k+2},\ldots ,
\nu_{n}$ for sections $H| _U$ can be chosen in such a manner that
brackets $[X,\nu_j]$ will be tangent vector fields to foliation
$F$ for each section $X$ of the bundle $H|_U$.Now assume that
foliation $F$ is a riemannian. Then for each tangent vector $ X\in
V(F)|_U$ it takes place $X_g(\nu_i,\nu_j)=0, \ i,j={k+1},\ldots
,n$ \cite{N1996}. By using this fact, it is easy to show that for
each vector field $ X\in V(F) $ takes place
\begin{center}
$Xg_{T}(Y,Z)=g_{T}([X,Y],Z)+g_{T}([Y,[X,Z]])$,
\end{center}
where $Y,Z\in V(M)$. In this case is $g_T$ called transversal
metrics  for foliation $F$, defined by the riemannian metric $g$
(\cite{MO1988}, see p.77). Notice that a transversal metrics $g_T$
determines the local distance between the leaves,since it defines
the length of the perpendicular geodesics.As follows from (
\cite{MO1988} the assertion 3.2), it is true also the  converse
fact i.e. if it is given $k$-dimensional foliation  $F$ on
riemannian manifold
 and riemannian metric $g$ defines transversal metric for $F$
then   $F$ is a riemannian foliation with  respect  to riemannian
metric $g$. Authors proved the similar fact for foliations with
singularities.

Let $F$ be a foliation with singularities, $L$ is a leaf of the
foliation $F$, $Q $ is a normal bundle of  $L$. Then riemannian
metric $g$ defines the metric $g^L_T$ on $Q$ as follows: if
$\nu_1,\nu_2:L\to Q$ are smooth (of the class $C^\infty$ )
sections of normal bundle $Q$ we will put
$g^L_T(\nu_1,\nu_2)=g(X,Y)$, where $X,Y\in V(M)$,the restrictions
of $X,Y$ on $L$ coincides with $\nu_1,\nu_2$ accordingly. The
metric $g^L_T$ is called transversal metric for $F$ on $L$, if for
each $X\in V(F)$ at
points of the leaf $L$ takes place\\
\begin{center}
$X{g^L_T}(Y,Z)=g^L_T([X,Y],Z)+g^L_T([Y,[X,Z]])$
\end{center}
where $Y,Z\in V(M)$, $g^L_T(Y,Z)=g(\pi Y,\pi Z)$ , $\pi:TM\to Q$
is the orthogonal projection  considered  over $L$. We will notice
that riemannian foliation with singularities has no
$n$-dimensional leaves.There is an assumption  in (\cite{MO1988},
p. 201) that if complete riemannian metric defines on each leaf
foliation $F$ transversal metric,then foliation $F$ will be a
riemannian. This problem is solved positively by following
theorem.

\textbf {Theorem-2.5.} \cite{N1996} Let $M$ be a complete
riemannian manifold with riemannian  metric $g$, $F$ is a singular
foliation on $M$ which  has no $n$- dimensional leaves. Then the
foliation $F$ is a riemannian if and only if riemannian metric $g$
defines on each leaf of the foliation $F$ transversal metric.

Now we will assume that $x_0\in M$, $L_0=L(x_0)$ is a proper leaf
with trivial holonomy group and the system (2.1) is completely
controllable on $L_0$.

\textbf {Theorem-2.6} \cite{AN2004} Let the mapping $x\to L(x)$ be
continuous at a point $x_0$. Then the system (2.1)is completely
controllable on the orbits,sufficiently close to $L_0$.

\textbf {Remark.} Multiple-valued mapping $x\to L(x)$ is called to
be lower semicontinuous at  a point $x_0$ if for each open set $V$
such that $V\bigcap L(x_0)\neq \emptyset $ there exists
neighborhood $B_{x_0}$  of the point $x_0$  such that $L(x)\bigcap
V\neq\emptyset$ for $x\in B_{x_0}$. Multiple-valued mapping $x\to
L(x)$ is called to be upper semicontinuous at a point $x_0$ if for
each open set  $V$ such that $L(x_0)\subset V $, there is a
neighborhood $B_{x_0}$  of the point $x_0$ such that $L(x)\subset
V$ for all $x\in B_{x_0}$. Multiple-valued mapping is continuous
at a point $x_0$ if it simultaneously lower and upper
semicontinuous at a point $x_0$.In our case is easy to prove that
mapping $x\to L(x)$ is  lower semicontinuous at each point of $M$
\cite{N1995}.

Sufficient conditions at which mapping $x\to L(x)$ is continuous,
are studied in papers of A.Narmanov \cite{N1995}, \cite{N1997}. We
will bring some of them. The following sufficient condition on
continuity of multiple-valued mapping $x\to L(x)$ follows directly
from the theorem 1.9 of the part one.

\textbf {Theorem-2.7}. Assume that $dimA_x(D)=k$ for every $x\in
M$, where $0<k<n$.If the set $L(x_0)$ is a compact leaf with
finite holonomy group then mapping $x\to L(x)$ is continuous at
the point $x_0$.

The following theorem shows that if foliation $F$ is a singular
riemannian foliation then the mapping $x\to L(x)$ is continuous at
each point $x_0$.

\textbf {Theorem-2.8} (\cite{N1995}) Let $F$ is a riemannian
foliation with singularities. Then multiple-valued mapping $x\to
L(x)$ is continuous at each point of the manifold $M$.

Now we will consider the problem on a continuity Bellman function
for a optimal time problem. We will remind that Bellman function
$T_\eta (x):G_\eta\to R^1$ is defined as follows: $T_\eta(\eta)
=0$ , $T_\eta(x)=inf{(\tau: there \ exists \ trajectory \
\alpha:[0,\tau ]\to M }$ ${of \ the \ system (2.1) \ such\ that \,
\alpha(0)=x, \ \alpha (\tau )=\eta)}$. The structure of set of
controllability generally can be rather difficult. Now we will
determine a class of control systems for which a set of
controllability  $G_\eta$ of the system (2.1) for all $\eta\in M$
coincides with a orbit $L(\eta)$ of family of the vector fields
$D=\{f(\cdot ,u):u\in U\}$ .

\textbf {Definition-2.4} \cite{NPI1985} We will say that the
system (2.1) is continuously-balanced at a point $x\in M$ if for
each vector field $X\in D$ there are vector fields $X_1,X_2,\ldots
,X_k$ from $D$, a neighborhood $V(x)$ of the point $x$ and the
positive continuous functions $\lambda_1(y),\lambda_2(y),\ldots
,\lambda_k(y)$ which are given in this neighborhood such that for
all $y\in V(x)$ takes place equality: $X(y)+\sum\
{\lambda_i(y)X_i(y)}=0.$

If we assume that the system (2.1) is continuously-balanced at
each point $ x \in M $ by means of results of work \cite{St1974}
it is possible to show that for each $\eta \in M $ the set of
controllability $ G_\eta $ of system  (2.1) coincides with the
orbit $L (\eta) $ of the family $ D=\{f(\cdot ,u):u\in U\} $ .

\textbf {Definition-2.5}. We will say that function $ T=T_\eta(x)$
is continuous at the point $ x_0\in G_\eta $ if for every
$\varepsilon>0 $ there is such neighborhood $V $of the point $ x_0
$ in topology of $ M $ that for any point $ x\in G_\eta\bigcap V $
takes place inequality  $\mid
T_\eta(x)-T_\eta(x_0)\mid<\varepsilon$.

For the system (2.1) given on compact manifold the following
result is obtained by professor N.N.Petrov \cite{NP1985}.

\textbf {Theorem-2.9}.Let $M$ be compact manifold. Then following
assertions are equivalent:

1) System  (2.1)  is $N$- locally controllable  near the  point
$\eta$.

2) For each $T>0$ the set $G_\eta(<T)$ is a domain in the manifold
$L(\eta)$.

3) For each $T>0$ the level    ${x\in M:T_\eta(x)=T}$ is the
border of the set$G_\eta(<T)$.

4) Bellman function $ T=T_\eta(x)$ is continuous at every point of
$G_\eta$.

Thus,for compact manifold the problem on a continuity of Bellman
function is reduced to a question about $N$ - local
controllability of system (2.1). For noncompact manifolds the
following theorem gives the  necessary and sufficient conditions
of a continuity of Bellman function  which is presented in
\cite{NP1985}.

\textbf {Theorem-2.10}. Let the system (2.1)is
continuous-balanced   at each point of $M$, for each $T>0$ the set
$G_\eta(\leq T)$ has compact closure and $dimA_x(D)=const$ for
every $x\in G_\eta$.Then  Bellman  function is continuous on
$G_\eta$ if and only if $G_\eta$ is a proper leaf of the foliation
generated by orbits of $D$.

In a case when manifold  $M$ is a analytic and the set $D$
consists of analytical vector fields owing to theorem Nagano
\cite{NAG1966} the condition $dimA_x(D)=const$ for all $x\in
G_\eta$ is always satisfied. Generally $dimA_x(D)$ can vary from a
point to a point on $G_\eta$ and always $ dimA_x(D)\leq dimL_\eta
$ for $x\in G_\eta$. For continuously-balanced control systems the
following theorem takes place \cite{N1997}.

\textbf {Theorem-2.11}. The set $G_\eta$ is a proper leaf of the
foliation generated by orbits of the family $D$ if and only if the
set $G_\eta$ is a set of type $F_\sigma$ and $G_\delta$
simultaneously in topology  of  manifold  $M$.

\bibliographystyle{plain}

\end{document}